\documentclass[11pt]{amsart}
\usepackage{amssymb}
 \newcommand{\Z}{{\mathbb Z}}
 \newcommand{\C}{{\mathbb C}}
 \newcommand{\F}{{\mathbb F}}
 \newcommand{\Q}{{\mathbb Q}}
 \newcommand{\R}{{\mathbb R}}
 \newcommand{\K}{{\mathbb K}}
\renewcommand{\L}{{\mathbb L}}
\renewcommand{\frak}{\mathfrak}
\newcommand{\ve}{\varepsilon}
 \newtheorem{theorem}{Theorem}[section]
 \newtheorem{lemma}[theorem]{Lemma} 
  \newtheorem{proposition}[theorem]{Proposition}
 \newtheorem{corollary}[theorem]{Corollary} 
  
 \newtheorem{example}[theorem]{Example}

\renewcommand{\o}{{\mathcal{O}}}
\newcommand{\wt}{\widetilde}
 \def\Box
  {\hfill \thinspace\vbox{\hrule height .5pt \hbox{\vrule  
   width .5pt \vbox to 7pt{\hbox to 3.5pt{}} \vrule width .5pt} 
   \hrule height 0pt depth .5pt}}

\title{Analogies between group actions on $3$-manifolds and number 
fields}

\author{Adam S. Sikora}



\begin{document}
\thispagestyle{empty}

\begin{abstract} Let a cyclic group $G$ act either on a number field $\L$
or on a $3$-manifold $M.$ Let $s_{\L}$ be the number of
ramified primes in the extension $\L^G\subset \L$ and $s_M$ be the number
of components of the branching set of the branched covering 
$M\to M/G.$ In this paper, we prove several formulas relating
$s_{\L}$ and $s_M$ to the induced $G$-action on $Cl(\L)$ and 
$H_1(M),$ respectively.
We observe that the formulas for $3$-manifolds and number fields are
almost identical, and therefore, they provide new evidence for
the correspondence between $3$-manifolds and number fields
postulated in arithmetic topology.

This is an extended version of a paper which will appear
in {\em Comm. Math. Helv.}
\end{abstract}

\maketitle

\tableofcontents

\section{Introduction and statement of the main results}
\pagestyle{myheadings} 
\markboth{\hfil{\sc ADAM S. SIKORA}\hfil}
{\hfil{\sc GROUP ACTIONS ON $3$-MANIFOLDS AND NUMBER FIELDS}\hfil}

Denote the ring of algebraic integers in a number field $\K$ 
by $\o_{\K}.$ Mazur's calculations of the \'etale cohomology
groups of rings of algebraic integers, \cite{Mazur-etale}, show that
for number fields $\K$ the groups $H^n_{et}(Spec\, \o_{\K},{\bf G}_m)$ 
vanish (up to $2$-torsion) for $n>3,$ and that they are equal to 
$\Q/\Z$ for $n=3.$
Since, furthermore, the groups $H^*_{et}(\o_{\K},{\bf G}_m)$ satisfy 
Artin-Verdier duality which is reminiscent of $3$-dimensional Poincar\'e 
duality, B. Mazur and D. Mumford suggested a
surprising analogy between the spaces $Spec\, \o_{\K}$ and
$3$-dimensional manifolds. The points in $Spec\, \o_{\K}$
(which are immersions of the finite residue fields $\o_{\K}/{\frak p}$
into $\o_{\K}$) can be viewed as $1$-dimensional objects and hence 
compared to knots in a $3$-manifold. Note that the fundamental group of 
a circle is $\Z$ and the absolute Galois group of a finite field is 
the profinite completion of $\Z.$
Further analogies between number fields and $3$-manifolds were
described by B. Mazur, A. Reznikov, M. Morishita and others 
in \cite{Mazur-60,Mo1,Mo2, Ramachandran, Reznikov1, Reznikov3, Waldspurger},
making a foundation for 
``arithmetic topology.'' 
At the heart of it lies a ``dictionary'' (which we call the
MKR dictionary after Mazur, Kapranov, and Reznikov)
matching the corresponding terms from $3$-dimensional topology and 
number theory, \cite{Reznikov2}. 
Despite its limitations and inconsistencies, the dictionary
can be used for translating statements from $3$-dimensional topology
into number theory, and vice versa, often with a surprising accuracy.
We describe the basic version of the MKR dictionary in Section 
\ref{sarithmetic_top}, compare \cite{Reznikov2}.
(Further details can be found in \cite{TCFT}.)

In this section we state the main results of this paper,
which were motivated by the following two problems:
\vspace*{.1in}

Let $C_p$ be a cyclic group of prime order $p.$\\

{\bf Problem T:} {\em  Let $C_p$ act on a smooth, closed, connected, oriented 
$3$-manifold $M$ by orientation preserving diffeomorphisms.
In this situation the projection $\pi: M\to M/C_p$ is a branched covering.
Given the induced $C_p$-action on the torsion and 
free parts of $H_1(M,\Z),$ denoted by $H_{tor}(M)$ and $H_{free}(M)=
H_1(M,\Z)/H_{tor}(M)$, find the best lower and upper estimates
on the number, $s,$ of components of the branching set.}

Upon translation into number theory, Problem T assumes the following form:

{\bf Problem N:} {\em Let $C_p$ act on a number field $\L$ and let
$\K=\L^{C_p}.$ Given the induced $C_p$-action on the ideal class group, 
$Cl(\L),$ and on the group of units, $\o_{\L}^*,$ 
find the best lower and upper estimates on the number of 
ramified primes in the extension $\K\subset \L.$
We will denote by $s_0$ and $s$ the numbers of finite ramified primes
and all ramified primes (including the infinite ones), 
respectively. (For $p\ne 2,$ $s_0=s$).}

Our particular interest in these two problems stems from the fact
that none of them has an elementary solution, and furthermore,
our solutions are based on methods which 
are beyond the current scope of arithmetic topology. Therefore, 
we hope that the consideration of these two problems will provide a new 
insight into arithmetic topology.

According to the MKR dictionary (Section \ref{dictionary}),
$H_{tor}(M)$ and $H_{free}(M)$ correspond to 
$Cl(\L)$ and $\o_{\L}^*/torsion,$ respectively.
For that reason, we will denote the latter two groups by $H_{tor}(\L)$ and 
$H_{free}(\L).$

Let $\F_p$ denote the field of $p$ elements. 

\begin{theorem}\label{upper}
(1) Under the assumptions of Problem T, and the additional assumption
$H_{free}(M/C_p)=0,$ we have
$$s\leq 1+dim_{\F_p}\,H^2(C_p, H_{tor}(M))+ 
dim_{\F_p}\,H^1(C_p,H_{free}(M)).$$
(2) Under the assumptions of Problem N 
$$s_0 \leq 1+ dim_{\F_p}\, H^2(C_p, H_{tor}(\L))+ 
dim_{\F_p}\, H^1(C_p, H_{free}(\L)).
$$
\end{theorem}

Recall that if $H$ is an abelian group with a $C_p$-action, 
$C_p=<\tau| \tau^p=1>,$ 
then $H^i(C_p,H)= Ker N/Im S$ for odd $i>0$ and $H^i(C_p,H)= Ker S/Im N,$ 
for even $i>0,$ where $N$ and $S$ are homomorphisms from $H$ to $H$ given 
by the multiplication by
$1+\tau+...+\tau^{p-1}$ and by $\tau-1,$ respectively.
Furthermore, the method of Herbrand quotient (\cite[Prop. 2.53]{Koch}), 
implies that
\begin{equation}
H^1(C_p, H)\simeq H^2(C_p, H) \label{H_q}
\end{equation}
for finite $H.$

We do not have a satisfying explanation for the coincidence of the formulas
of Theorem \ref{upper}. Their proofs (given in Sections \ref{s_cyclic_3mfld} 
and \ref{s_cyclic_nf}) are based on two very different methods.
Examples \ref{eg_Hempel} and \ref{eg_lens} show that the extra 
assumption of Theorem \ref{upper}(1) is necessary and its inequality cannot
be improved. Similarly, the inequality of (2) cannot be improved -- see
comments following Theorem \ref{lower}.

Here is an upper bound for $s,$ which does not need the additional assumption
of Theorem \ref{upper}(1).

\begin{theorem}\label{upperT}
Under the assumptions of Problem T,
$$s\leq 1+dim_{\F_p}\,H^2(C_p, H_1(M))+ 
dim_{\F_p}\,H^1(C_p,H_{free}(M)).$$
\end{theorem}

Searching for lower estimates for $s$ 
one encounters the following problem:
A $G$-action on a $3$-manifold $M$ induces $G$-actions on
$H_1(M),$ $H_{tor}(M),$ and $H_{free}(M).$ However, $H_1(M,\Z)$ and 
$H_{free}(M)\oplus H_{tor}(M)$ do not need to be isomorphic as $G$-modules!
(See Example \ref{eg_Hempel}).
Since this problem does not occur in number theory, 
the next result requires that
$H_1(M,\Z)= H_{free}(M)\oplus H_{tor}(M)$ as $G$-modules.

\begin{theorem} \label{lower}
(1)  If conditions of Problem T are satisfied,
$s>0,$ and \\ $H_1(M,\Z)=H_{free}(M)\oplus H_{tor}(M)$ as 
$C_p$-modules then
$$s\geq 1+dim_{\F_p} H^2(C_p,H_{tor}(M)).$$
(2) If conditions of Problem N are satisfied,
$s> 0,$ and $Cl(\K)$ has no elements of order $p$ then
$$s\geq 1+dim_{\F_p} H^2(C_p,H_{tor}(\L)).$$
\end{theorem}

\begin{corollary} \label{corlower}
(1)  If $H_{tor}(M/C_p)$ has no  elements of order $p$
then  $$\text{$p$-tor} \left(H_{tor}(M)^{C_p}\right)=(\Z/p)^d,$$ 
for some $d.$
($p$-tor denotes the $p$-torsion part of the group.)
Furthermore, under the conditions of Theorem \ref{lower}(1),
$d\leq s-1.$

(2) If $H_{tor}(\K)$ has no elements of order $p$ then
$$\text{$p$-tor} \left(H_{tor}(\L)^{C_p}\right)=(\Z/p)^d.$$
Furthermore, under the conditions of Theorem \ref{lower}(2),
$d\leq s-1.$
\end{corollary}

\begin{proof}(1) Note that 
$$H^2(C_p,H_{tor}(M))=\text{$p$-tor}\left(H_{tor}(M)^{C_p}\right)/ 
\text{$p$-tor}\left(Im\, N\right).$$
Now the result follows from the fact that
$N=\pi_\sharp \pi_*:H_{tor}(M)\to H_{tor}(M),$ 
where $\pi_*:H_{tor}(M)\to H_{tor}(M/C_p),$ cf. Sec. \ref{ss_transfer}.\\
The proof of (2) is identical.
\end{proof}

To illustrate these results we consider quadratic number fields.
Let $\L=\Q(\sqrt{d})$ where $d>0$ is a square-free integer.
$Gal(\L/\Q)=C_2$ acts on $Cl(\L)$ by the involution $I\to I^{-1}.$
Hence, by (\ref{H_q}), $H^1(C_2, Cl(\L))=H^2(C_2, Cl(\L))=Cl(\L)^{C_2}.$
A simple calculation yields $H^1(C_2,H_{free}(\L))=\F_2$ and, 
therefore, Theorems \ref{upper}(2) and \ref{lower}(2) imply the following
estimates (similar to a formula of Gauss):
$$1\leq s-dim_{\F_2}\, Cl(\L)^{C_2}\leq 2.$$
These estimates cannot be improved. In fact $s-dim_{\F_2}\, Cl(\L)^{C_2}$
is either $1$ or $2$ depending if the norm of fundamental unit in $\o_{\L}^*$
is $-1$ or $1$ ($d=10$ represents the first case, $d=15$ the second).

The topological examples are more laborious and, therefore, they are 
postponed to Section \ref{sexamples}. Example \ref{eg_Hempel} shows
that the extra assumption in Theorem \ref{lower}(1) 
(saying that $H_1(M,\Z)=H_{free}(M)\oplus H_{tor}(M)$ as $C_p$-modules) 
is necessary and that the inequality of Theorem \ref{lower}(1) is sharp.

The situation with cyclic extensions of $\Q$ of dimensions greater than $2$
is much simpler. The following generalization of the Gauss formula
follows from \cite[Lemma 13.4.1]{Lang-cyclo}; see also our proof enclosed in
Section \ref{s_cyclic_nf}.

\begin{theorem}\label{generalized_Gauss}
If $\L$ is a Galois extension of $\Q$ of prime degree $p\ne 2$ then
$$H_{tor}(\L)^{C_p}=(\Z/p)^{s-1}.$$ 
\end{theorem}


The topological counterpart of the above formula is the following restatement
of a result due to Reznikov, \cite[Thm. 15.2.5]{Reznikov3}.
This result follows from our Theorems \ref{upper}(1) and \ref{lower}(1) 
and from the fact that $H_{free}(M)=0$ implies $H_{free}(M/C_p)=0$ by
inequality (\ref{b1}) in Section \ref{ss_transfer}.

\begin{corollary}\label{reznikov-thm}
If $M$ is a rational homology sphere, $H_{tor}(M/C_p)=0,$ 
and $s\ne 0$ then
$$H_{tor}(M)^{C_p}=(\Z/p)^{s-1}.$$ 
\end{corollary}

Example \ref{eg_lens} shows that the assumption $H_{free}(M)=0$ is necessary.

Despite the fact that all above formulas are pairwise identical
we do not know any uniform proof of them. In fact, the statements 
concerning group actions on manifolds are proved
by calculations on the Leray-Serre spectral sequence for 
equivariant cohomology, and the statements concerning
Problem N are proved using local methods (ideles) and class 
field theory. This suggests that the results of class field theory used in 
the proof can be interpreted in terms of $3$-dimensional topology. 
(This suggestion was formulated before, for example by B. Mazur and 
A. Reznikov).
Analogously, one might expect that equivariant cohomology
can be formulated and used in the framework of number theory.

{\bf Acknowledgments:} We would like to thank C. Allday, J. Hempel,
N. Ramachandran, and L. Washington for helpful conversations.
Additionally, we express our special thanks to R. Schoof
for suggesting the method of proving Theorems \ref{upper}(2) and 
\ref{lower}(2), and the referee of this paper for pointing out \cite{CH}
as a basis of an alternative approach to proving these results.

\section{Arithmetic topology}
\label{sarithmetic_top}

\subsection{Algebraic number fields}
\label{snumberth}

Let $\K$ be an algebraic number field, i.e. a field which is 
a finite extension of $\Q.$
The elements of $\K$ which are roots of monic polynomials in $\Z[x]$ are
called algebraic integers and form a ring denoted by $\o_{\K}.$
For example, $\o_{\Q}=\Z,$ $\o_{\Q(\sqrt{d})}=\Z[\sqrt{d}]$ 
for $d\ne 1$ mod $4,$ and $\o_{\Q(\sqrt{d})}=\Z[\frac{1+\sqrt{d}}{2}]$ 
otherwise.
The two most important groups associated with $\K$ are the group of
units of $\o_{\K}$ denoted by $\o_{\K}^*$ and the ideal class
group, $Cl(\K).$ The latter group is composed of all non-zero 
ideals in $\o_{\K}$ considered up to an equivalence 
$\sim,$ such that $I\sim J$ if and only if $I=c\cdot J,$ for some 
$c\in \K^*.$
The multiplication of ideals in $\o_{\K}$ provides $Cl(\K)$ with a 
structure of an abelian group. This group is always finite.
$\o_{\K}$ is a principal ideal domain if and only if
$Cl(\K)$ is trivial. Since principal ideal domains coincide with
unique factorization domains in the class of Dedekind domains 
(which includes the rings $\o_{\K}$),
$Cl(\K)$ measures the extend to which these two properties fail for $\o_{\K}.$

An introduction to algebraic number theory can be found for example
in \cite{FT, Ja}.

\subsection{MKR dictionary}
\label{dictionary}
\label{skap-rez}
Kapranov and Reznikov, building on ideas of Mazur, proposed the 
following dictionary matching corresponding terms from $3$-dimensional
topology and number theory:
\begin{enumerate}
\item Closed, oriented, connected, smooth $3$-manifolds correspond to 
algebraic number fields.
(More precisely, $3$-manifolds correspond to the schemes $Spec\, \o_{\K},$
where $\K$'s are algebraic number fields).
\item A link in $M$ corresponds to an ideal in $\o_{\K}$ and a
knot in $M$ corresponds to a prime ideal
in $\o_{\K}.$ (We consider only tame links and knots).
Note that, knots can be represented by immersions
of $S^1$ into $M$, and
prime ideals in $\o_{\K}$ can be identified with closed immersions
$Spec\, \F\to Spec\,\o_{\K},$ where $\F$'s are finite fields.
Each link decomposes uniquely
as a union of knots and each ideal decomposes uniquely as a product of primes.
\item An algebraic integer $w\in \o_{\K}$ is analogous to an embedded surface
(possibly with boundary) $S\subset M.$ The operation $w\to (w)\triangleleft 
\o_{\K}$ corresponds to $S\to \partial S.$
Note that ideals of the form $(w)$ represent the identity in 
$Cl(\K),$ and the links of the 
form $\partial S$ represent the identity in $H_1(M,\Z)$.
\item A closed embedded surface $F\subset M$ is analogous to a unit
in $\o_{\K}.$
\item $Cl(\K)$ corresponds to the torsion part of $H_1(M,\Z),$ denoted by
$H_{tor}(M).$
\item the group of units in $\o_{\K}$ modulo the torsion (the roots
of unity), $\o_{\K}^*/\mu(\K),$ corresponds to $H_{free}(M)=
H_1(M,\Z)/H_{tor}(M).$
\item Finite extensions of number fields $\K\subset \L$ correspond to
finite branched coverings of $3$-manifolds $\pi:M\to N.$
(The branching set may be empty.)
\item Each Galois extension $\K\subset \L$ with 
$Gal(\L/\K)=G$ induces the morphism
$Spec\,\o_\L\to (Spec\,\o_\L)/G= Spec\, \o_\K.$ Such maps correspond to 
the quotient maps $M\to M/G$ induced by orientation preserving
actions of finite groups $G$ on $3$-manifolds $M.$ 
One can show that $M/G$ is always a $3$-manifold and that the maps 
$M\to M/G$ are branched coverings. (However, if $G$ is not cyclic then
one needs to consider here a broad class of branched coverings 
whose sets of ramified points may be graphs.)
\item $S^3$ corresponds to $\Q.$ 
The $3$-sphere has no non-trivial (unbranched) covers. Similarly, 
$\Q$ has no non-trivial unramified extensions. (Ramified extensions
are defined in Section \ref{ss_split_ramified_inert}.)
\item Let $\xi_{p^n}$ be a primitive $p^n$-th root of $1,$ where $p$ is
prime. The extension $\Q\subset \Q(\xi_{p^n})$ is ramified at $p$ only.
Such extensions correspond to the cyclic branched covers of knots in $S^3.$ 
The infinite cyclotomic extension $\Q\subset \Q(\xi_{p^\infty})=
\bigcup_n \Q(\xi_{p^n})$ corresponds
to the infinite cyclic cover of a knot complement, 
$\wt{S^3\setminus {\mathcal K}}.$ 
There is a natural action of $\Z$ on $H_1(\wt{S^3\setminus {\mathcal K}})$ 
and a natural action of the $p$-adic integers on the $p$-torsion of
$Cl(\Q(\xi_{p^\infty}))=\lim_{\leftarrow} Cl(\Q(\xi_{p^n})).$  These actions
give rise to the Alexander polynomial of $\mathcal K$ and 
the characteristic polynomial of $p,$ cf. \cite{Mazur-60}.
For an exposition of Iwasawa theory including more information on 
characteristic polynomials see \cite{Wa}.
\end{enumerate}

This basic version of the dictionary is based on the preprint 
\cite{Reznikov2}. (In its published version, \cite{Reznikov3}, 
Reznikov introduces a different dictionary in which number fields 
correspond to ``$3{1\over 2}$-manifolds.'' We will not discuss it here
since our results do not support his new dictionary.)
Our more precise version of the dictionary can be found in \cite{TCFT}.

Despite the above analogies, one should not expect to find a bijection 
between number fields and $3$-manifolds. Hence, $\Q(\sqrt{5})$ does not 
correspond to any specific $3$-manifold. Furthermore, we will see
later that the existence of a bijection between number fields and 
$3$-manifolds fulfilling the conditions of the MKR dictionary would 
disprove Poincar\'e conjecture, cf. \cite{Ramachandran}.
Nevertheless, we believe that there exists a theory 
which encompasses (at least partially) both algebraic number theory and 
$3$-dimensional topology.

\subsection{Galois extensions and Galois branched covers, I}\ \\
\label{ss_GaloisI}
Let $G=Gal(\L/\K)$ and let ${\frak q}$ be a prime ideal in $\o_{\L}.$ 
The decomposition group of $\frak q,$
$D_{\frak q}\subset G,$ is
the group of elements of $G$ preserving ${\frak q},$
$$D_{\frak q}=\{g\in G: g({\frak q})={\frak q}\}.$$
The quotient $\o_{\L}/{\frak q}$ is a finite field and the image of the 
homomorphism $D_{\frak q}\to Gal(\o_{\L}/{\frak q})$ is 
composed of exactly those automorphisms of $\o_{\L}/{\frak q}$ 
which fix the subfield 
$\o_{\K}/{\frak p},$ where ${\frak p}=\o_{\K}\cap {\frak q}.$
The kernel of this homomorphism, $I_{\frak q},$ is called the inertia 
group of ${\frak q}.$
Hence, we have the following exact sequence,
$$0\to I_{\frak q}\to D_{\frak q}\to  Gal\left(\o_{\L}/{\frak q}/
\o_{\K}/{\frak p}\right)\to 0.$$
The order of $I_{\frak q},$ denoted by $e_{\frak q}$ is 
called the ramification index of ${\frak q}.$
The order of $Gal\left(\o_{\L}/{\frak q}/
\o_{\K}/{\frak p}\right)$ will be denoted
by $f_{\frak q}.$  

Assume now that a finite group $G$ acts on a $3$-manifold $M$ and 
that $\pi:M\to M/G$ is a branched covering (ie. assume that the branching 
set of $\pi$ is a link).
The subgroup  $D_{\mathcal K}\subset G$ composed of elements mapping a 
knot ${\mathcal K}\subset M$ to itself will be called the 
decomposition group of $\mathcal K.$ Assume that
the action of $D_{\mathcal K}$ on ${\mathcal K}$ is orientation preserving.
In this situation the image of the natural homomorphism $D_{\mathcal K}\to 
{\rm Homeo}({\mathcal K})$ is exactly the group of deck transformations, 
$Gal({\mathcal K}/{\mathcal K'}),$ of the covering 
${\mathcal K}\to {\mathcal K}'={\mathcal K}/D_{\mathcal K}.$ 
By analogy with the number theory, the kernel of this homomorphism,
$I_{\mathcal K},$  is called the inertia group of
$\mathcal K.$ As before, we have the following exact sequence:
$$0\to I_{\mathcal K}\to D_{\mathcal K}\to 
Gal({\mathcal K}/{\mathcal K'})\to 0.$$
We denote the orders of $I_{\mathcal K}$ and of 
$Gal({\mathcal K}/{\mathcal K'})$ 
by $e_{\mathcal K}$ and $f_{\mathcal K}.$ 
We will call $e_{\frak q}$ the ramification index of ${\mathcal K}.$
Note that both $Gal\left(\o_{\L}/{\frak q}/
\o_{\K}/{\frak p}\right)$ and $Gal({\mathcal K}/{\mathcal K'})$
are always cyclic groups.

\subsection{Galois extensions and Galois branched covers, II}\ \\
\label{ss_GaloisII}
As before, let $\frak q$ be a prime ideal in $\o_{\L}$ and let 
${\frak p}={\frak q}\cap \o_{\K},$ $G=Gal(\L/\K).$
The ideal ${\frak p}\o_{\L}\triangleleft \o_{\L}$
decomposes uniquely as a product of prime ideals, 
\mbox{${\frak p_1}^{e_1}\cdot ... \cdot {\frak p_g}^{e_g},$}
where $e_i$ is the 
ramification index of $\frak p_i$ (defined in the previous section) 
and $\frak q$ is one of the ideals ${\frak p_1},...,{\frak p_g}.$

\begin{theorem}(see eg. \cite{FT} Theorem III.20) \label{field_ext}
Under the above assumptions, 
\begin{itemize}
\item $G$ acts transitively on ${\frak p_1}, ...,{\frak p_g},$
\item $e_1=...=e_g\stackrel{def}{=} e$ and $f_{{\frak p}_1}=... =
f_{{\frak p}_g}\stackrel{def}{=}f;$
\item $|G|=efg.$
\end{itemize}
\end{theorem}

Let $G$ act on a $3$-manifold $M$ such that $\pi: M\to M/G$ is a 
branched covering. Consider a knot $\mathcal K$ in $M$ such that 
${\mathcal K}/G$ is a knot in
$M/G$ which is either a component of the branching set or it is disjoint
from the branching set.
In this situation, $\pi^{-1}(\mathcal K)$ is a link in $M$ whose components
we denote by ${\mathcal K}_1, ..., {\mathcal K}_g,$
$$\pi^{-1}(\mathcal K)={\mathcal K}_1\cup ... \cup {\mathcal K}_g.$$

We leave to the reader the proof of the following theorem analogous to
Theorem \ref{field_ext}

\begin{theorem} Under the above assumptions,
\begin{itemize}
\item $G$ acts transitively on ${\mathcal K}_1, ..., {\mathcal K}_g;$
\item $e_{{\mathcal K}_1}=...=e_{{\mathcal K}_g} \stackrel{def}{=} e;$ 
$f_{{\mathcal K}_1}=... =f_{{\mathcal K}_g}\stackrel{def}{=}f;$
\item $|G|=efg.$
\end{itemize}
\end{theorem}

\subsection{Split, ramified, and inert primes and knots.}
\label{ss_split_ramified_inert}
A prime $\frak q\triangleleft \o_{\L}$ is ramified in a Galois extension
$\K=\L^G\subset \L$ if $e_{\frak q}>1,$
$\frak q$ is split if $e_{\frak q}=f_{\frak q}=1,$ and $\frak q$ is
inert if $e_{\frak q}=1, f_{\frak q}=|G|.$

Consider now a Galois branched covering $\pi:M\to M/G$ and a knot
${\mathcal K}\subset M$ satisfying the assumptions of the previous section.
By analogy to the number theoretic terminology,
we will say that $\mathcal K$ is ramified in the branched covering 
$\pi$ if $e_{\mathcal K}>1,$ $\mathcal K$ splits if 
$e_{\mathcal K}=f_{\mathcal K}=1,$ and $\mathcal K$ is inert if 
$e_{\mathcal K}=1, f_{\mathcal K}=|G|.$ Observe that ${\mathcal K}$ is
ramified if and only if it is a component of the branching set.

If $G=C_p$ then each prime ${\mathfrak q}\triangleleft \o_{\L}$ and each 
knot ${\mathcal K}\subset M$ is either split, inert, or ramified. If 
${\mathfrak p}= {\mathfrak q}\cap \o_{\K}$ then ${\mathfrak q}$ is
\begin{itemize}
\item split if ${\frak p} \o_\L={\frak q}_1\cdot ...\cdot {\frak q}_p,$
where ${\frak q}_1, ..., {\frak q}_p,$ are different primes, one of which is
${\frak p}.$ In this situation $C_p$ cyclically permutes these primes.
\item ramified if ${\frak p} \o_\L={\frak q}^p.$ 
In this situation $C_p$ fixes the elements of ${\frak q}.$
\item inert if ${\frak p} \o_L={\frak q}.$
In this situation $C_p$ acts non-trivially on ${\frak q}.$
\end{itemize}

If $G=C_p$ then $\pi: M\to M/G$ is a branched covering. If 
${\mathcal K}\subset M$ satisfies the assumptions of the previous section
then ${\mathcal K}$ is
\begin{itemize}
\item split if $\pi^{-1}({\mathcal K}/G)={\mathcal K}_1\cup ...
\cup {\mathcal K}_p,$ where ${\mathcal K}_1, ..., {\mathcal K}_p$
are different knots one of which is $\mathcal K.$
In this situation $C_p$ cyclicly permutes these knots.
\item ramified if ${\mathcal K}/G$ is a component of a branching set.
In this situation $C_p$ fixes $\pi^{-1}({\mathcal K}/G)={\mathcal K}.$
\item inert if $\pi^{-1}({\mathcal K}/G)={\mathcal K}$ and the $C_p$-action
on $\mathcal K$ is non-trivial. ($C_p$ rotates $\mathcal K$ by the angle
$2\pi/p$).
\end{itemize}

\begin{theorem}\label{similar_properties} 
Consider a Galois extension $\K\subset \L,$ $G=Gal(\L/\K),$ and
a Galois branched covering $\pi: M\to M/G.$
\begin{enumerate}
\item There are only finitely many ramified primes in $\o_{\L}$ and ramified 
knots in $M;$
\item There are infinitely many split primes in $\o_{\L}$ and infinitely many 
split knots in $M;$
\item If $G$ is cyclic then there are infinitely many inert primes in 
$\o_{\L}$ and infinitely many inert knots in $M.$
\item If $G$ is not cyclic then there are no inert primes nor inert knots.
\item If $G$ is cyclic of prime order then each knot and prime is either 
split, ramified, or inert. (Recall that we consider only knots 
${\mathcal K}\subset M$ such that ${\mathcal K}/G$ is a knot
in $M/G$ which is either a component of the branching set or it is
disjoint from the branching set).
\end{enumerate}
\end{theorem}

The proof of the topological statements of Theorem \ref{similar_properties}
is left to the reader.
{\it Proof of number theoretic statements:}\\
(1) follows from \cite[Theorem 7.3]{Ja}.\\
(2) ${\frak q}$ is split if and only if $D_{\frak q}=\{1\},$ 
and hence, if and only if, $\frak q$ is not 
ramified 
and its Frobenius, $\left[\frac{\L/\K}{\frak q}\right],$ is trivial.
By Chebotarev's density theorem, \cite[Thm. 1.116]{Koch}, there are 
infinitely
many such primes.\\
(3) If $G$ is cyclic then by Chebotarev's density theorem there are infinitely
many unramified primes ${\frak q}\triangleleft \o_{\L}$ whose 
Frobenius is a generator of $G$ and, hence, $D_{\frak q}=G.$ 
Each such ${\frak q}$ is inert.\\
(4) If $\frak q$ is inert then $G=D_{\frak q}=Gal(\L/{\frak q}/\K/{\frak p})$
is cyclic.\\
(5) By Theorem \ref{field_ext}, $efg=p.$ Hence
exactly one of the indexes $e,f,g$ is equal to $p.$
\Box

Note that the MKR dictionary together with 
the analogy between split, ramified, and inert primes and knots is
sufficient to translate a version of Poincar\'e conjecture
into a statement in number theory. The Poincar\'e conjecture
states that $S^3$ is the only closed $3$-manifold with no unbranched covers.
However, the geometrization conjecture implies a stronger statement: 
$S^3$ is the only closed $3$-manifold with no finite unbranched covers. 
Somewhat surprisingly, its algebraic 
translation: ``$\Q$ is the only number field with no unramified extensions''
is false. However, only a few counterexamples are known and it is not
known whether there are only a finite number of them.

\begin{theorem}(see eg. \cite[Ex. 11.2]{Wa}) If $d<0$ and 
the class group of $\L=\Q(\sqrt{d})$ is trivial then $\L$ has no unramified 
extensions. There are precisely
$9$ such values of $d:$ $-1,-2,-3,-7,-11,-19,-43,-67,-163.$
\end{theorem}

\subsection{Transfer and norm maps}
\label{ss_transfer}
For any $G$-action on $M$ consider the map $\pi_*:H_1(M,\Z)\to H_1(M/G,\Z)$
induced by $\pi:M\to M/G$ and the transfer map
$\pi_\sharp: H_1(M/G,\Z)\to H_1(M,\Z)$ defined as follows: 
For any $x\in H_1(M/G,\Z)$ represented by a closed curve $c$ disjoint 
from the branching
set, $\pi_\sharp(x)=[\pi^{-1}(c)]\in H_1(M,\Z).$

The induced maps
$$\pi_*:H_{free}(M,\Z)\to H_{free}(M/G,\Z),\ 
\pi_*:H_{tor}(M,\Z)\to H_{tor}(M/G,\Z),$$
$$\pi_\sharp: H_{free}(M/G)\to H_{free}(M),\  
\pi_\sharp: H_{tor}(M/G)\to H_{tor}(M),$$
satisfy the following conditions:
\begin{enumerate}
\item $\pi_*\pi_\sharp$ is the multiplication by $|G|.$
\item $\pi_\sharp\pi_*(x)=\sum_{g\in G} g\cdot x.$
\end{enumerate}
In particular we have
\begin{equation}\label{b1}
rank\, H_{free}(M/C_p)\leq rank\, H_{free}(M).
\end{equation}

Given a Galois extension $\K\subset \L,$ with $Gal(\L/\K)=G,$
we have the norm map $\pi_*: \o_{\L}\to \o_{\K},$ 
$\pi_*(x)=\prod_{g\in G} gx$ and the
embedding $\pi_\sharp:\o_{\K}\hookrightarrow \o_{\L}.$
The induced maps $\pi_*: H_{free}(\L)\to H_{free}(\K)$ and
$\pi_\sharp: H_{free}(\K)\to H_{free}(\L)$ satisfy properties
identical to those above:
\begin{enumerate}
\item $\pi_*\pi_\sharp(x)=x^{|G|}.$
\item $\pi_\sharp\pi_*(x)=\prod_{g\in G} gx$
\end{enumerate}
Furthermore, we can define $\pi_\sharp$ and $\pi_*$ for ideals:
$\pi_\sharp(I)=I\cdot \o_{\L}$ for $I\triangleleft \o_{\K},$ and
$\pi_*(J)=\prod_{g\in G} gJ$ for  $J\triangleleft \o_{\L}.$
The induced maps $$\pi_\sharp: Cl(\K)\to Cl(\L),\quad
\pi_*: Cl(\L)\to Cl(\K),$$ once again satisfy the conditions
\begin{enumerate}
\item $\pi_*\pi_\sharp(I)=I^{|G|},$ 
\item $\pi_\sharp\pi_*(I)=\prod_{g\in G} gI.$
\end{enumerate}

\section{Cyclic group actions on manifolds}
\label{s_cyclic_3mfld}
In this section we will prove the results concerning cyclic group
actions on $3$-manifolds which were announced in the introduction.
As references for the homological algebra used in this section 
(group cohomology and spectral sequences) we recommend 
\cite{Brown, Weibel}.

\subsection{Algebraic preliminaries}
Let $\Z_{(p)}$ denote the ring of integers localized at $(p),$
$p$ prime.

\begin{proposition}\label{classification}
If $V$ is a free abelian group and an indecomposable $\Z[C_p]$-module
then $V\otimes \Z_{(p)}$ is isomorphic to one of the following 
$\Z_{(p)}[C_p]$-modules:
\begin{enumerate}
\item[(F)] the free module $\Z_{(p)}[C_p],$
\item[(T)] the trivial module, $\Z_{(p)},$ or
\item[(AI)] the augmentation ideal, $Ker\, \ve,$ for $\ve:\Z_{(p)}[C_p] \to
\Z_{(p)},$ $\ve(g)=1.$
\end{enumerate}
\end{proposition}

The proof follows from 
the classification of indecomposable $C_p$-representations over 
$\Z$, due to Diederichsen and Reiner, \cite[Theorem \S 74.3]{Curtis}.\Box

For any $C_p$-module $V$ the cohomology groups $H^i(C_p,V)$ for $i>0$ depend
on the parity of $i$ only and 
$$H^{2n}(C_p, V)=\hat{H}^{0}(C_p, V),\quad H^{2n-1}(C_p, V)=
\hat{H}^{1}(C_p, V),\quad \text{for\ } n>0,$$
are linear spaces over the field of $p$-elements, $\F_p.$
Here $\hat{H}^*(C_p, V)$ denotes the Tate cohomology groups.

The type of a $\Z[C_p]$-module $V$
(as described in Proposition \ref{classification}) can be determined
by the Tate cohomology of $C_p$ with coefficients in $V.$
If $V$ is a free abelian group and an indecomposable $\Z[C_p]$-module
then $V$ is of type 
\begin{enumerate}
\item[(F)] iff $\hat{H}^{i}(C_p, V)=0$ for $i=0,1,$
\item[(T)] iff $\hat{H}^0(C_p, V)=\Z/p,$ $\hat{H}^1(C_p, V)=0,$
\item[(AI)] iff $\hat{H}^0(C_p, V)=0,$ $\hat{H}^1(C_p, V)=\Z/p.$
\end{enumerate}

Any $C_p$-module $V$ gives rise to two other $C_p$-modules, $V^*$ and
$V^{\#},$ defined as follows: $V^*=Hom(V,\Z),$ and 
$g\in C_p$ acts on $Hom(V,\Z)$ by sending $f(\cdot)$ to $f(g^{-1}\cdot).$
The second module, $V^{\#},$ is equal to $V$ as 
an abelian group and the action of $g\in C_p$ on $v\in V^{\#}$ is given by 
$g^{-1}v.$ 

\begin{lemma} \label{dual_module} For any $C_p$-module $V,$\\
(1) $\hat{H}^*(C_p,V^{\#})=\hat{H}^*(C_p,V)$\\
(2) if $V$ is a free abelian group then $\hat{H}^*(C_p,V^*)=\hat{H}^*(C_p,V).$
\end{lemma}

\begin{proof} (1) As before, let $\tau$ be a generator of $C_p$
and let $N,S: V\to V$ be given by $\cdot \sum_{i=0}^{p-1} \tau^i$ and
by $\cdot(\tau-1)$ respectively. $\hat{H}^n(C_p,V)$ is equal to 
$Ker\, N/Im\, S$ or $Ker\, S/Im\, N,$
depending if $n$ is odd or even.
Similarly, $\hat{H}^n(C_p,V^{\#})$ is equal to 
$Ker\, N^{\#}/Im\, S^{\#}$ or $Ker\, S^{\#}/Im\, N^{\#},$
depending if $n$ is odd or even, where
$N^{\#},S^{\#}: V\to V$ are given by $\cdot \sum_{i=0}^{p-1} \tau^{-i}$ and
by $\cdot(\tau^{-1}-1)$ respectively. 
Now (1) follows from the fact that $N^{\#}=N,$ $Ker\, S^{\#}= Ker\, S,$ 
and $Im\, S^{\#}=Im\, S.$

(2) Let $V=\Z^d.$ Since $\hat{H}^*(C_p,V)$ is a $p$-group, $Tor(
\hat{H}^*(C_p,V),\Z_{(p)})=0$ (see \cite[Calculation 3.1.1]{Weibel}) and
by the universal coefficient theorem,
$$\hat{H}^i(C_p,V\otimes \Z_{(p)})=
\hat{H}^i(C_p,V)\otimes \Z_{(p)}=\hat{H}^i(C_p,V).$$
Therefore, it is sufficient to show that $V\otimes \Z_{(p)}$ and 
$Hom(V,\Z)\otimes \Z_{(p)}$ are isomorphic as $C_p$-modules, and
for that it is enough to consider indecomposable $C_p$-modules $V$ only.
Such modules are classified in Proposition \ref{classification}.
We leave it to the reader to check that in each of the three possible 
cases we get 
$$V\otimes\Z_{(p)}\simeq Hom(V\otimes\Z_{(p)},\Z_{(p)})
\simeq Hom(V,\Z)\otimes\Z_{(p)}$$ as needed.
\end{proof}

\subsection{Equivariant cohomology and Leray-Serre spectral sequence}

Let $(X,A)$ be a relative $C_p$--CW-complex. In other words, let $X$ be a 
CW-complex with a $C_p$-action, $A$ be a  subcomplex, and let the
$C_p$-action on $X$ preserve
$A$ and be compatible with the CW-structure on $X.$
Additionally, we require that if $g\in C_p$
maps a cell $\sigma\subset X$ into itself 
then $g$ acts trivially on $\sigma;$ compare \cite{Allday-Puppe}. 
For any smooth $C_p$-action on a manifold $M,$
$(M,M^{C_p})$ is a relative $C_p$--CW-complex, \cite{Illman}.

Denote the cochain complex for the cellular cohomology of $(X,A)$ 
with coefficients in $\Z$ by $(C^*(X,A), \delta).$ 
Let $(P^*, \delta')$ be a complete projective resolution of
$\Z$ considered as a $\Z[C_p]$-module 
with the trivial $C_p$-action, cf. \cite[Ch. VI.3]{Brown}.
Then $D^{**}(X,A)=Hom_{\Z[C_p]}(P^*,C^*(X,A))$ is a double complex
with two differentials (a) $\delta_v=\delta: D^{kl}(X,A)\to D^{k,l+1}(X,A)$ 
and (b) $\delta_h: D^{kl}(X,A)\to D^{k+1,l}(X,A)$ dual to $\delta'.$
The cohomology of the associated total complex,
$$D^s=\bigoplus\limits_{k+l=s} D^{kl},$$
$$d(\alpha)=d_h(\alpha)+(-1)^kd_v(\alpha),\ \text{for $\alpha\in D^{kl}$}$$ 
is the Tate equivariant cohomology.

Consider the ``first'' spectral sequence $(E^{**}_*(X,A),d_*)$ associated with 
$(D^{**}(X,A),\delta_h,\delta_v),$ induced by the vertical filtration,
cf. \cite{Swan-new}.
Its first term is
$E^{kl}_1(X,A)=H^l(X,A;\Z),$ for all $k,l\in \Z.$

The next two results describing the properties of $(E^{**}_*(X,A),d_*)$
belong to the folk knowledge.

\begin{lemma} \label{0-row} If 
$X$ is a connected $C_p$--CW-complex and $X^{C_p}\ne \emptyset,$ 
then all differentials $d_r^{*,r-1}: E^{*,r-1}_r(X)\to 
E^{*0}_r(X)$ are $0$ for $r\geq 2.$ 
(In other words, the $0$th row of $E^{**}_2(X)$ 
survives to infinity).
\end{lemma}

\begin{proof}
Let $\{*\}$ denote the one point space and let $x_0$ be a zero cell in 
$X^{C_p}.$ The $C_p$-equivariant maps: $*\to x_0: \{*\}\to X$ and 
$X\to \{*\}$ induce maps
$$(E^{**}_r(*), d_r)\to (E^{**}_r(X), d_r)\to 
(E^{**}_r(*), d_r),$$
whose composition is the identity.
$E_2^{i0}(X)=H^i(C_p,H^0(X,\Z))$ is either $0$ or $\F_p.$
Hence, if $d_r^{k,r-1}: E^{k,r-1}_r(X)\to E^{k+r,0}_r(X)$ is not $0$ for
some $k,r$ then $E^{k+r,0}_{r+1}(X)=0$ and
$E^{k+r,0}_r(X)=\F_p.$ Hence, $k+r$ is even and, hence, 
$E^{k+r,0}_{r+1}(*)=\F_p.$
Therefore $$E^{**}_{r+1}(*)\to E^{**}_{r+1}(X)\to E^{**}_{r+1}(*)$$
cannot be the identity.
\end{proof}

\begin{lemma} \label{n-row} If $(M,B)$ is a $C_p$--CW-complex such that
$M$ is a connected, closed, oriented $n$-manifold and $B$ is 
an $n$-dimensional ball then all differentials 
$d_r^{*,n}: E^{*,n}_r(M)\to E^{*,n-r+1}_r(M)$
are $0$ for $r\ge 2.$ (In other words, the $n$-th row of $E^{**}_2(M)$ 
survives to infinity).
\end{lemma}

\noindent{\it Proof (due to T. Skjelbred).} 
The sequence of cochain complexes
$$C^*(M,M\setminus int\, B)\to C^*(M)\to C^*(M\setminus int\, B)$$
induced by the embedding $M\setminus int\, B\hookrightarrow M$ 
yields a sequence of spectral sequences 
$$E^{**}_r(M,M\setminus int\, B)\stackrel{\alpha_r}{\to} E^{**}_r(M)\to 
E^{**}_r(M\setminus int\, B)$$
which is exact for $r=1.$
Since $E^{*n}_1(M\setminus int\, B)=H^n(M\setminus int\, B)=0,$
the map $\alpha^{*n}_1: E^{*n}_1(M,M\setminus int\, B)\to E^{*n}_1(M)$ 
is onto.

Since $$E^{*k}_1(M,M\setminus int\, B)=H^k(M,M\setminus int\, B)=
H^k(B^n,\partial B^n)=
H_{n-k}(B^n),$$
all non zero elements of $E^{**}_r(M,M\setminus int\, B)$ for $r\geq 1$ 
lie in the $n$th row. Since $\alpha_r$
commutes with differentials and all differentials in 
$E^{**}_r(M,M\setminus int\, B)$ are $0$ for $r\geq 2,$
$d_r^{*,n}: E^{*,n}_r(M)\to E^{*,n-r+1}_r(M)$ is zero for $r\geq 2.$
\Box

\subsection{Cyclic group actions on $3$-manifolds}
Let $M$ be a smooth, connected, closed, oriented
$3$-manifold with a smooth orientation preserving $C_p$-action with
a fixed point.

By the result of Illman \cite{Illman} recalled above,
$M$ can be given a structure of a relative 
$C_p$--CW-complex. Furthermore, since $M^{C_p}\ne \emptyset,$
one may assume that there exists a $C_p$--CW-subcomplex $B$ of $M$ 
homeomorphic to a $3$-ball.
In this situation, the second term of $E^{**}_*(M)$ is:
\[
\begin{array}{ccccc}
...& \F_p & 0 & \F_p & ...\\
... & \hat{H}^0(C_p,H^2(M))&
\hat{H}^1(C_p,H^2(M))& \hat{H}^0(C_p,H^2(M)) & ...\\

... & \hat{H}^0(C_p,H^1(M))&
\hat{H}^1(C_p,H^1(M))& \hat{H}^0(C_p,H^1(M)) & ...\\

...& \F_p & 0 & \F_p & ...\\
\end{array}\]

By Poincar\'e duality $H^2(M,\Z)=H_1(M,\Z)^{\#},$ $H^1(M,\Z)=H_{free}(M)^*$
as $C_p$-modules, and by Lemma \ref{dual_module},
\begin{equation}
\label{dual}
\begin{array}{c}
\hat{H}^n(C_p, H^2(M,\Z))\simeq \hat{H}^n(C_p, H_1(M,\Z)),\\
\hat{H}^n(C_p, H^1(M,\Z))\simeq \hat{H}^n(C_p, H_{free}(M,\Z)).
\end{array}
\end{equation} 

Let 
$$\Psi^n: \hat{H}^n(C_p,H_1(M,\Z))\to \hat{H}^n(C_p,H_{free}(M))$$
denote the homomorphism given by the composition of 
$d_2^{n2}: E^{n,2}_2\to E^{n+2,1}_2$ with the isomorphisms (\ref{dual}).
One needs to be aware that $\Psi^n$ is not the map induced by the natural
homomorphism $H_1(M,\Z)\to H_{free}(M).$

\begin{theorem}\label{s} If $M^{C_p}\ne \emptyset$ ($p$ prime) then
$M^{C_p}$ is a union of $s$ circles, where
$$s=1+dim_{\F_p}\, Ker\, \Psi^n+dim_{\F_p}\, 
Coim\, \Psi^{n-1},$$
for any $n.$
\end{theorem}

Before giving the proof we note that this result immediately implies 
Theorem \ref{upperT}.
It also implies Theorems \ref{upper}(1) and \ref{lower}(1)
\vspace*{.1in}:

{\bf Proof of Theorem \ref{upper}(1):} If $H_{free}(M/C_p)=0$ then by
property (2) of the transfer map (cf. Subsection \ref{ss_transfer})
$H_{free}(M)$ is annihilated by $\sum_{g\in C_p} g.$
By Proposition \ref{classification},
$H_{free}(M)\otimes \Z_{(p)}$ is a sum of three types of indecomposable 
$\Z_{(p)}[C_p]$-modules. Since $Ker\, \ve$ is the only indecomposable 
$\Z_{(p)}[C_p]$-module which is annihilated by $\sum_{g\in C_p} g,$ 
$H_{free}(M)\otimes \Z_{(p)}$ is a sum of
modules of this type. We compute that $H^2(C_p, Ker\, \ve)=0.$
Since $\otimes \Z_{(p)}$ is an exact functor in the
category of $C_p$-modules, $H^2(C_p, H_{free}(M))=
H^2(C_p, H_{free}(M)\otimes \Z_{(p)})=0.$ Now, the short exact sequence 
$$0\to H_{tor}(M)\to H_1(M)\to H_{free}(M)\to 0$$
implies that $$dim_{\F_p}\, H^2(C_p,H_1(M))\leq dim_{\F_p}\, 
H^2(C_p,H_{tor}(M)).$$ Therefore, Theorem \ref{upper}(1)
follows from Theorem \ref{upperT}.
\vspace*{.1in}
\Box

{\bf Proof of Theorem \ref{lower}(1):}
If $H_1(M,\Z)=H_{free}(M)\oplus H_{tor}(M)$ as $C_p$-modules,
the dimension of $Ker\, \Psi^n$ is bounded below by
$$\begin{array}{c}
dim_{\F_p} \hat{H}^n(C_p,H_1(M))- dim_{\F_p} \hat{H}^n(C_p,H_{free}(M))\\
=dim_{\F_p} \hat{H}^n(C_p,H_{tor}(M)).
\end{array}$$\Box

We complete this section by proving Theorem \ref{s}.

By Lemmas \ref{0-row} and \ref{n-row} all differentials in the
spectral sequence $(E^{*,*}_r(M),d_r),$ for $r\geq 2,$ are
$0$ except possibly $d_2^{k2}: E^{k2}_2(M)\to E^{k+2,1}_2(M).$ 
Therefore, $$\sum_{k+l=n} E^{kl}_\infty(M)=
\F_p+ Ker\, d_2^{n-2,2}+ Coim\, d_2^{n-3,2}\simeq\\
\F_p+ Ker\, \Psi^{n-2}+ Coim\, \Psi^{n-3}.$$ 
Since $(E^{**}_*(M),d_*)$ converges to $\hat{H}_{C_p}^*(M,\Z),$
$\hat{H}_{C_p}^n(M,\Z)$ and $\sum_{k+l=n} E^{kl}_\infty(M)$ have equal numbers
of elements. By the Localization Theorem, \cite[Ch. VII Prop. 10.1]{Brown},
$\hat{H}_{C_p}^n(M,\Z)=\hat{H}_{C_p}^n(M^{C_p},\Z).$
Since $M^{C_p}$ is composed of circles,
Theorem \ref{s} follows from the following lemma.

\begin{lemma} $\hat{H}^n_{C_p}(M^{C_p},\Z)=\F_p^s,$ for any $n.$
\end{lemma}

\begin{proof}
Since the $C_p$-action on $M^{C_p}$ is trivial,
the total complex of $(D^{**}(M^{C_p}),\delta,\delta')$ is isomorphic to
$(P^*)\otimes (C^*(M^{C_p},\Z),\delta^*)$ where $(P^*)$
is the cochain complex 
$\stackrel{0}{\to} \Z\stackrel{\cdot p}{\to}\Z\stackrel{0}{\to} \Z
\stackrel{\cdot p}{\to} ...$ Hence, by  K\"unneth formula,
$\hat{H}^n_{C_p}(M^{C_p},\Z)= H^n(D^{**}(M^{C_p}))$
is equal to 
$$\bigoplus_{k+l=n}
H^k(P^*)\otimes H^l(M^{C_p},\Z)+
\bigoplus_{k+l=n+1} Tor_{\Z}(H^k(P^*), H^l(M^{C_p},\Z)).$$
Since $H^k(P^*)$ is either $\F_p$ or $0$ depending if $k$ even or odd, the
proof follows.
\end{proof}

\section{Cyclic group actions on number fields}
\label{s_cyclic_nf}
In this section we prove Theorems \ref{upper}(2),
\ref{lower}(2), and \ref{generalized_Gauss} stated in the
introduction. The proofs use class field theory and they are based on 
the approach explained to us by R. Schoof. Although most of the necessary
definitions will be stated in the course of the proof, we refer the
reader unfamiliar with class field theory to \cite{CF,Koch,Ja,Lang-ant} for
more background information.

Let $I(\K)$ denote the group of fractional ideals in $\K$. Its elements are
finitely generated $\o_{\K}$-submodules of $\K,$ and the group product of 
$I_1$ and $I_2$ is $I_1\cdot I_2\subset \K.$
The fractional ideals generated by a single element are called principal and
they form a subgroup $P(\K)\subset I(\K).$ 

A valuation $\nu$ on $\K$ is a homomorphism $\nu: \K^* \to  \R^*_{>0},$ 
such that $\nu(x+y)\leq \nu(x)+\nu(y)$ for $x,y\in \K^*, x+y\ne 0.$
The valuation $\nu=1$  is called trivial, and
from now on it will be ignored. Valuations $\nu_1,\nu_2$
such that $\nu_2(x)=\nu_1(x)^c,$ for some $c>0$ are called equivalent and
the equivalence classes of valuations are called places.

Any prime $\mathfrak p$ in $\o_{\K}$ induces the $\mathfrak p$-adic
valuation, $\nu_{\mathfrak p}:\K^*\to (0,\infty)$ such that
$\nu_{\mathfrak p}(x)=d^{-\upsilon_{\mathfrak p}(x)},$
where $d$ is any number greater than $1$ and $\upsilon_{\mathfrak p}(x)$ is 
the power of $\mathfrak p$ in the prime decomposition of the fractional 
ideal $(x).$ (Usually, $d=|\o_{\L}/{\mathfrak p}|.$)
Such valuations are called non-archimedean.
Another type of valuations is induced by embeddings of $\K$ into $\R$ or 
$\C.$ If $\imath$ is such an embedding then $\nu(x)=|\imath(x)|$
is a valuation on $\K$ called archimedean.

One studies properties of fields by considering either prime ideals in 
$\o_{\K}$ or valuations on $\K.$ 
These two approaches are almost equivalent, since every valuation
is given either by a prime or by an embedding, as above.
For that reason, the embeddings $\imath: \K\hookrightarrow \R$
and the pairs of conjugate embeddings $\imath, \bar{\imath}: 
\K\hookrightarrow \C$ are called infinite primes. 
The theory of extensions of valuations 
is analogous to the theory of extensions of prime ideals.

Let $C_p$ act on a field $\L$ and let $\K=\L^{C_p}.$
Let $s_0$ denote the number of finite ramified primes (see Sec. 
\ref{ss_split_ramified_inert}) in the extension $\K\subset \L$ and let 
$s_{\infty}$ 
denote the number of infinite ramified primes. An infinite prime 
$\imath, \bar{\imath}: \L\hookrightarrow \C$ is ramified if
$\imath(\K) \subset \R$ and $\imath(\L)\not \subset \R.$
Infinite ramified primes may appear for $p=2$ only. Let $s$ be the 
number of all ramified primes, $s=s_0+s_{\infty}.$

Consider the group $\prod_{\nu} \L^*_{\nu},$ where the product is taken 
over all places of $\L$ and $\L_{\nu}^*$ is the group of non-zero
elements in the completion of $\L$ with respect to $\nu.$
The group of ideles, ${\mathbb I}_{\L},$ is the subgroup of 
$\prod_{\nu} \L^*_{\nu}$ composed of the elements $(x_{\nu})$ such that 
$x_{\nu}\in (\o_{\L})^*_{\nu}$ for almost all places $\nu.$ 

Consider the embedding $\L^*\to {\mathbb I}_{\L}$ mapping 
$x$ to $(x_{\nu})$ where $x_{\nu}$ is the image of $x$ under the 
completion map $\L\to \L_{\nu}.$ The quotient, 
$C(\L)={\mathbb I}_{\L}/\L^*,$ is called the idele class group of $\L.$

Observe that if $\nu$ is non-archimedean then 
${\mathfrak p}_{\nu}=\{x\in \K^* | \nu(x)<1\}$ is a
prime ideal in $\o_{\K}$ and $\nu$ is equivalent to the 
${\mathfrak p}_{\nu}$-adic
valuation. The completion of ${\mathfrak p}_{\nu}$ in 
$(\o_{\L})_{\nu}$ is the only prime ideal in $(\o_{\L})_{\nu}.$ 
Hence for each $x\in \L^*_{\nu}$ there exists $\upsilon_{\nu}(x)\in \Z$
such that $(x)={\mathfrak p}_{\nu}^{\upsilon_{\nu}(x)}.$ 
Note that there is a homomorphism ${\mathbb I}_{\L}\to I(\L)$ 
sending $(x_{\nu})$ to 
$\prod_{\nu} {\mathfrak p}_{\nu}^{\upsilon_{\nu}(x_{\nu})},$
where the product is taken over all non-archimedean places.
The kernel of this homomorphism, $U_\L,$ is called the group of 
idelic units. Notice that $\o_\L^*$ embeds into $U_\L.$

We have the following exact and commuting diagram in which
all maps except the bottom horizontal ones are either 
defined above or obvious.

\begin{equation}
\begin{array}{ccccccc}\label{diagram}
      & 0         &     & 0          &    & 0         & \\
      & \downarrow&     & \downarrow &    &\downarrow & \\
0 \to & \o_\L^*    & \to & \L^*        & \to& P(\L)      & \to 0\\
      & \downarrow&     & \downarrow &    &\downarrow & \\
0 \to & U_\L       & \to &{\mathbb I}_\L&\to& I(\L)      & \to 0\\
      & \downarrow&     & \downarrow &    &\downarrow & \\
0 \to & U_\L/\o_\L^*& \to & C(\L)       & \to& Cl(\L)     & \to 0\\
      & \downarrow&     & \downarrow &    &\downarrow & \\
      & 0         &     & 0          &    & 0         & \\
\end{array}
\end{equation}

The existence of the bottom horizontal maps follows from the 
following lemma.  The proof of it is left to the reader.
\vspace*{.1in}

\begin{lemma}
If the following diagram is exact and commutes
\[
\begin{array}{ccccccc}
      & 0         &     & 0          &    & 0         & \\
      & \downarrow&     & \downarrow &    &\downarrow & \\
0 \to & A_{02}    & \to & A_{12}     & \to& A_{22}    & \to 0\\
      &\downarrow &     & \downarrow &    &\downarrow& \\
0 \to & A_{01}    & \to & A_{11}     & \to& A_{21}    & \to 0\\
      &\downarrow &     & \downarrow &    &\downarrow & \\
      & A_{00}    &     & A_{10}     &    & A_{20} & \\
      &\downarrow &     & \downarrow &    &\downarrow & \\
      & 0         &     & 0          &    & 0         & \\
\end{array}\]
then there exist unique bottom horizontal maps making the diagram
an exact commuting $3\times 3$ diagram.
\end{lemma}

Note that $C_p$ acts  on each of the $9$ groups in
(\ref{diagram}).
We need to calculate the cohomology groups of $C_p$ with 
coefficients in these groups.

\begin{lemma}\label{U_L}
$\hat{H}^i(C_p,U_\L)=\begin{cases}\F_p^s & \text{for $i=0$}\\ 
\F_p^{s_0} & \text{for $i=1.$}\\
\end{cases}$
\end{lemma}

\begin{proof}
$$U_\L=\prod_{\nu\ arch} \hspace*{-.1in} \L_{\nu}^*\cdot 
\prod_{\nu\ nonarch} \o_{\L,\nu}^*$$
is equal as a $C_p$-module to
$$\prod_{\nu'} M_{\nu'}$$
where the product is taken over all places $\nu'$ in $\K$ and
$$M_{\nu'}=\begin{cases}
\prod_{\nu |\nu'} \o_{\L,\nu}^* & \text{if $\nu'$ is non-archimedean}\\
\prod_{\nu |\nu'} \L_{\nu}^* & \text{if $\nu'$ is archimedean.}\\
\end{cases}$$
Since $C_p$ is a simple group, a place $\nu$ in $\K$ either splits
completely in $\L,$ or it is inert, or it is totally ramified.
In the first case $C_p$ acts on $M_{\nu'}$ by freely permuting its 
components and therefore $$\hat{H}^*(C_p,M_{\nu'})=0.$$
In the second and third cases $M_{\nu'}$ is either $\o_{\L,\nu}^*$
or $\L_{\nu}^*$ and the cohomology groups of $C_p$ with these coefficients 
have been calculated in Lemma 4 in \cite[Ch. IX \S3]{Lang-ant}:
$\hat{H}^*(C_p,M_{\nu'})=0,$ if $\nu'$ is inert, and
$$\hat{H}^0(C_p,M_{\nu'})=\F_p,\quad \hat{H}^1(C_p,M_{\nu'})=
\begin{cases} 0 & \text{if $\nu'$ is archimedean}\\ 
\F_p & \text{if $\nu'$ is non-archimedean}\\
\end{cases}$$
for ramified $\nu'.$ (Note that for archimedean valuations,
the calculation reduces to calculating $\hat{H}^*(C_2,\C)$ where
$C_2$ acts on $\C$ by conjugation).
Hence $$\hat{H}^i(C_p,U_\L)=\prod_{\nu'}
\hat{H}^i(C_p,M_{\nu'})=\F_p^s \ {\rm or\ } \F_p^{s_0},$$
depending if $i=0$ or $i=1.$
\end{proof}

One of the main results of class field theory is the calculation
of the cohomology groups of $C_p$ with coefficients in the idele
class group (see \cite[Ch. IX \S 5]{Lang-ant}): 
$$\hat{H}^i(C_p,C(\L))=\begin{cases}
\F_p & \text{for $i$ even}\\ 
0 & \text{for $i$ odd.}\\
\end{cases}$$
The isomorphism $\hat{H}^0(C_p,C(\L))=C(\K)/NC(\L)\to C_p=\F_p$ is given
by the Artin map. If $\omega$ ramifies in the extension $\K\subset \L,$ 
$\nu=\omega^p,$ then by the corollary after Theorem 4 \cite[Ch. XI]{Lang-ant} 
the Artin map restricted to $\hat{H}^0(C_p,\o_{\L,\omega}^*)$
is also an isomorphism
onto $C_p.$ Therefore, if $s>0$ then
the composition of the maps $\o_{\L,\omega}^*\to U_\L\to {\mathbb I}_\L\to 
C(\L),$ induces an epimorphism $$\hat{H}^0(C_p,\o_{\L,\omega}^*)\to 
\hat{H}^0(C_p,C(\L)).$$
Hence, from diagram (\ref{diagram}) we get:

\begin{corollary}\label{epi}
If $s>0$ then the map $$\hat{H}^0(C_p,U_\L/\o_\L^*)\to \hat{H}^0(C_p,C(\L))$$
is an epimorphism.
\end{corollary}

Therefore the bottom row of (\ref{diagram}) yields the exact sequence
$$\hat{H}^0(C_p,U_\L/\o_\L^*) \twoheadrightarrow  \hat{H}^0(C_p,C(\L))
\stackrel{0}{\to} \hat{H}^0(C_p,Cl(\L))\to \hat{H}^1(C_p,U_\L/\o_\L^*)\to 0$$
and, hence, $$\hat{H}^0(C_p,Cl(\L))=\hat{H}^1(C_p,U_\L/\o_\L^*).$$
On the other hand, the left column of (\ref{diagram}) gives the exact sequence
$$\hat{H}^1(C_p,\o_\L^*)\to \hat{H}^1(C_p,U_\L)=\F_p^{s_0}\to
\hat{H}^1(C_p,U_\L/\o_\L^*)$$ and, therefore,
\begin{equation}\label{best_estimate}
s_0\leq dim_{\F_p}\,\hat{H}^0(C_p,Cl(\L))+ dim_{\F_p}\,\hat{H}^1
(C_p,\o_\L^*).
\end{equation}
Let $\mu_*$ denote the torsion part of $\o_\L^*$ ($\mu_*$ is the group of 
roots of $1$ in $\L$). Since $\mu_*$ is cyclic, 
$dim_{\F_p}\,\hat{H}^1(C_p,\mu_*)\leq 1.$ 
Finally, the sequence $$\hat{H}^1(C_p,\mu_*)\to
\hat{H}^1(C_p,\o_\L^*)\to \hat{H}^1(C_p,\o_\L^*/\mu_*)$$
is exact, and therefore we have
$$dim_{\F_p}\,\hat{H}^1(C_p,\o_\L^*)\leq 1+ dim_{\F_p}\,
\hat{H}^1(C_p,H_{free}(\L)).$$
(Recall that $H_{free}(\L)=\o_\L^*/\mu_*$).
This inequality together with (\ref{best_estimate}) completes the proof of 
Theorem \ref{upper}(2).

In order to prove Theorem \ref{lower}(2) we need another fact from class 
field theory:

\begin{proposition}\label{shoof} (due to R. Schoof)
$$Ker(Cl(\K)\stackrel{\pi_{\sharp}}{\rightarrow} Cl(\L))=
Ker(H^1(C_p,\o_\L^*)\to H^1(C_p,U_\L)).$$
($\pi_{\sharp}: Cl(\K)\to Cl(\L)$ was defined in Sec. \ref{ss_transfer}).
\end{proposition}

\begin{proof}
Since the functor $N\to N^{C_p}$ on $C_p$-modules is left exact,
the bottom horizontal row of (\ref{diagram}) yields the exact sequence
$$0 \to  (U_{\L}/O_{\L}^*)^{C_p} \to C(\L)^{C_p}  \to Cl(\L)^{C_p}.$$
By considering a $3\times 3$ diagram for $\K$ analogous to (\ref{diagram})
we get the exact sequence
$$0\to U_{\K}/\o_{\K}^*\to C(\K)\to Cl(\K)\to 0$$ which together with
the sequence above forms a commutative diagram
\begin{equation}\label{snake}
\begin{array}{ccccccc}
0 \to & U_\K/\o_\K^*&      \to & C(\K)       & \to& Cl(\K)     & \to 0\\
      & \downarrow i&       & \downarrow &    &\downarrow \pi_\sharp& \\
0 \to & (U_\L/O_\L^*)^{C_p}& \to & C(\L)^{C_p}&  \to & Cl(\L)^{C_p}.& \\
\end{array}
\end{equation}
By Proposition 8.1 in \cite{CF}, the central vertical arrow is an isomorphism.
Hence, by Snake Lemma (cf. \cite[Lemma 1.3.2]{Weibel}), 
\begin{equation}\label{ker=coker}
Ker\, \pi_\sharp=Coker\, i.
\end{equation}
Consider the long exact sequence of cohomology groups of $C_p$ induced by the
short exact sequence of coefficient groups
$$0\to \o_\L^* \to U_\L \to U_\L/ \o_\L^*\to 0.$$ 
Since for any $C_p$-module $N,$ $H^0(C_p,N)=N^{C_p},$ and $(\o_\L^*)^{C_p}=
\o_\K^*,$ $(U_\L)^{C_p}=U_\K,$
 this long exact sequence gives rise to an exact sequence
$$0 \to U_\K/O_\K^* \stackrel{i}{\to}  (U_\L/O_\L^*)^{C_p} 
\stackrel{g}{\to} H^1(C_p, O_\L^*) \stackrel{h} \to H^1(C_p, U_\L)\to ...$$
where $i$ coincides with the map $i$ of (\ref{snake}).
Hence, by (\ref{ker=coker}), 
$$Ker\, \pi_\sharp=Coker\, i=(U_\L/O_\L^*)^{C_p}/Ker\ g= Im\, g=Ker\, h.$$
\end{proof}

{\bf Proof of Theorem \ref{lower}(2):} 
If $Cl(\K)$ has no elements of order $p$ then, by Proposition \ref{shoof},
$H^1(C_p,\o_\L^*)$ embeds into $H^1(C_p,U_\L)$ and 
hence the map $\hat{H}^0(C_p,U_\L/\o_\L^*)\to H^1(C_p,\o_\L^*)$ is $0.$
Therefore, $\hat{H}^0(C_p,U_\L)\to \hat{H}^0(C_p,U_\L/\o_\L^*)$ is an 
epimorphism and by Lemma \ref{U_L}
\begin{equation}\label{nf-1} 
s\geq dim_{\F_p}\,\hat{H}^0(C_p,U_\L/\o_\L^*).
\end{equation}
On the other hand, the bottom row of (\ref{diagram}) and Corollary \ref{epi}
yield the exact sequence
$$0=\hat{H}^1(C_p,C(\L)) \to \hat{H}^1(C_p,Cl(\L))\to 
\hat{H}^0(C_p,U_\L/\o_\L^*) \twoheadrightarrow  \hat{H}^0(C_p,C(\L))=\F_p.$$
Hence 
\begin{equation}\label{nf-2}
dim_{\F_p}\, \hat{H}^0(C_p,U_\L/\o_\L^*)=1+ \,
dim_{\F_p} \hat{H}^1(C_p,Cl(\L)).
\end{equation}
Finally, by (\ref{H_q}), $\hat{H}^1(C_p,Cl(\L))=\hat{H}^0(C_p,Cl(\L)).$
Therefore (\ref{nf-1}) and (\ref{nf-2}) imply
Theorem \ref{lower}(2).
\Box

{\bf Proof of Theorem \ref{generalized_Gauss}:}
It is enough to prove that $H^1(C_p, \o_{\L}^*)=\F_p,$ since then
the statement follows from (\ref{best_estimate}) and Theorem \ref{lower}(2).
Since the rank of $H_{free}(\L)$ is at most $p-1,$ 
Proposition \ref{classification} implies that 
$H_{free}(\L)\otimes \Z_{(p)}=Ker\, \ve.$
Therefore, $H^1(C_p,H_{free}(\L))=H^1(C_p,H_{free}(\L)\otimes \Z_{(p)})=\F_p.$
Additionally, since $\mu_{\L}^{C_p}=\{\pm 1\},$ we have 
$\hat{H}^0(C_p,\mu_{\L})=0,$ and, by (\ref{H_q}), 
$\hat{H}^1(C_p,\mu_{\L})=0.$ Now it follows
from the exact sequence $0\to \mu_{\L}\to \o_{\L}^*\to H_{free}(\L)\to 0$ 
that $H^1(C_p, \o_{\L}^*)=\F_p.$
\Box

\section{Examples}
\label{sexamples}
We end this paper with two examples of group actions on $3$-manifolds
illustrating various anomalies related to cyclic group actions on 
$3$-manifolds.
The first example shows that the inequality of Theorem \ref{upper}(1) is sharp
and that the assumption $H_{free}(M)=0$ of Corollary \ref{reznikov-thm}
is necessary.

\begin{example}\label{eg_lens}
Consider a $C_p$-action on the lens space $L(p,1)$ with two circles
of fixed points. ($L(p,1)$ is obtained from two handlebodies identified
along their boundaries. $C_p$ acts on each of these solid tori by 
rotation along their cores. These actions give rise to the 
desired $C_p$-action on $L(p,1).$)
Let $U$ be an open
ball in $L(p,1)$ such that $U\cap gU=\emptyset$ for all $g\ne e$ in $C_p.$

Consider now a $C_p$-action on $S^3$ with one circle of fixed points.
Let $V$ be an open ball in $S^3$ such that $V\cap gV=\emptyset$ for all
$g\ne e$ in $C_p,$ and let $\Psi: \partial V\to \partial U$ be any 
homeomorphism.
By removing the balls $gU$ from $L(p,1)$ and the balls $gV$ from $S^3$
and by identifying the boundary spheres of $L(p,1)\setminus \bigcup_g gU$
and $S^3\setminus \bigcup_g gV$ by the homeomorphisms
$g\Psi g^{-1}: g\partial V\to g\partial U,$ we obtain a
closed $3$-manifold $M$ with a $C_p$-action.
Since $L(p,1)/C_p=S^3$ and $S^3/C_p=S^3,$ we have $M/C_p=L(p,1)/C_p\#
S^3/C_p=S^3.$
Simple calculations show that $H_1(M)=\F_p\oplus 
\Z[\xi]/(1+\xi+...+\xi^{p-1}),$ where $C_p$ acts on 
$\Z[\xi]/(1+\xi+...+\xi^{p-1})$ by $\cdot \xi$ and it acts trivially on
$\F_p.$
Hence $H^2(C_p,H_{tor}(M))=\F_p$ and $H^1(C_p,H_{free}(M))=\F_p.$
In this example $s=3.$
\end{example}

The second example shows that the assumption $H_{free}(M/C_p)=0$
in Theorem \ref{upper}(1) is necessary. Furthermore, it shows
that the assumption in Theorem \ref{lower}(1) saying that 
$H_1(M)=H_{tor}(M)\oplus H_{free}(M)$ as $C_p$-modules is also necessary.
(Take $n=1,$ $k=p$ below).

\begin{example}\label{eg_Hempel} (based on an idea of J. Hempel)\\
Let $\phi$ be a Dehn twist on the torus $T$ and let $M_k$ be the 
$3$-manifold obtained from $T\times [0,1]$ by identifying 
$T\times \{0\}$ with $T\times \{1\}$ via $\phi^k.$ Assume that
$\phi$ fixes neighborhoods of points $p_1,...,p_n$ and consider knots
$\{p_i\}\times S^1$ with the framing $v_i\times S^1,$ where $v_i$ is an
arbitrary tangent vector to $T$ at $p_i.$
By performing the surgery along these framed knots we obtain a new manifold 
$M_{k,n}.$
If $p$ divides $k$ then the map $(x,t)\to 
(\phi(x)^{k/p},t+1/p\ \text{mod}\ 1)$ defines a $C_p$-action on
$M_k$ which survives the surgery.
We have $H_1(M_k,\Z)=\Z\oplus \Z \oplus \Z/k$ and 
$H_1(M_{k,n},\Z)=\Z^n\oplus \Z/k,$ for $n\geq 1.$ Furthermore, 
$H_1(M_{p,n},\Z)$ decomposes as a $C_p$-module 
into a sum $\Z^{n-1}\oplus \Z\oplus \Z/p,$ with the trivial action
on the first summand and the action $$\Z\oplus \Z/p \to \Z\oplus \Z/p,\quad
(x,y)\to (x,x+y)$$ on the second summand. For this action
$H^2(C_p,H_{tor}(M_{p,n}))=\F_p,$\\
$H^1(C_p,H_{free}(M_{p,n}))=0,$ 
and $s=n.$ Note also that $M_{p,n}/C_p=M_{1,n}.$
\end{example}

\centerline{Dept. of Mathematics, 244 Math. Bldg.}
\centerline{University at Buffalo, SUNY}
\centerline{Buffalo, NY 14260, USA}
\centerline{asikora@buffalo.edu}


\end{document}